\documentclass[10pt,psamsfonts]{amsart}

\newtheorem{thm}{Theorem}[section]
\newtheorem{rk}[thm]{Remark}
\newtheorem{prop}[thm]{Proposition}
\newtheorem{clly}[thm]{Corollary}
\newtheorem{lemma}[thm]{Lemma}

\newtheorem{defi}[thm]{Definition}

\title[A curvature form for pseudoconnections]{A curvature form for pseudoconnections}

\author[C.A. Morales and M. Vilches]{C.A. Morales and M. Vilches}



\address{
C.A. Morales\\
Instituto de Matematica\\ Universidade Federal do Rio de Janeiro, 
P. O. Box 68530, 21945-970 Rio de Janeiro, Brazil.}
\address{
M. Vilches\\
Departamento de An\'alise Matem\'atica\\
IME, Universidade do Estado do Rio de Janeiro, 20550-013,
Rio de Janeiro, Brazil}

\begin{document}

\maketitle

\begin{abstract}
We obtain the curvature form $F^\nabla=P\circ d^\nabla\circ\nabla-d^\nabla\circ P\circ\nabla+d^\nabla\circ\nabla\circ P$
for a vector bundle pseudoconnection $\nabla$, where $d^\nabla$ is the exterior derivative associated to
$\nabla$.
We use $F^\nabla$ to obtain the curvature of $\nabla$.
We also prove that
$F^\nabla=0$ is a necessary (but not sufficient)
condition for $d^\nabla$ to be a chain complex.
Instead we prove that $F^\nabla=0$ and $d^\nabla\circ d^\nabla\circ\nabla=0$
are necessary and sufficient conditions for $d^\nabla$ to be
a {\em chain $2$-complex}, i.e., $d^\nabla\circ d^\nabla\circ d^\nabla=0$.
\end{abstract}

\section{Introduction}

\noindent
Let $M$ be a differentiable manifold.
Denote by $\Omega^0(M)$
the ring of $C^\infty$ real valued maps in $M$.
Denote by $\chi^\infty(M)$ and $\Omega^k(M)$
respectively the $\Omega^0(M)$-modules of 
$C^\infty$ vector fields and $k$-forms defined on $M$,
$k\geq 0$.
Denote by $d:\Omega^k(M)\to \Omega^{k+1}(M)$ the standard exterior derivation of $k$-forms of $M$.
We denote by $Hom(A,B)$ the set of homomorphisms from the modulus
$A$ to the modulus $B$.
If $\xi$ is a real smooth vector bundle over $M$ we denote by
$\Omega^k(\xi)$ the set of $\xi$-valued $k$-forms of $M$, namely,
$\Omega^0(\xi)$ is the $\Omega^0(M)$-module of
smooth sections of $\xi$ and
$\Omega^k(\xi)=\Omega^k(M)\otimes \Omega^0(\xi)$ for all $k\geq 1$.
If $\eta$ is another vector bundle over $M$ we denote by
$HOM(\xi,\eta)$ the set of bundle homomorphims
from $\xi$ to $\eta$ over the identity.
Every $P\in HOM(\xi,\eta)$ induces a homomorphism
$P\in Hom(\Omega^0(\xi),\Omega^0(\eta))$
of $\Omega^0(M)$ modules in the usual way.
It also defines a homomorphism $P\in Hom(\Omega^k(\xi),\Omega^k(\eta))$ for all $k\geq 1$ by
setting $P(\omega\otimes s)=\omega\otimes P(s)$ at every
generator $\omega\otimes s\in\Omega^k(\xi)$.

A {\em pseudoconnection} of a vector bundle $\xi$ over $M$
is an $I\!\!R$-linear map $\nabla:\Omega^0(\xi)\to\Omega^0(\xi)$
for which there is a bundle homomorphism $P\in HOM(\xi,\xi)$ called the principal homomorphism
of $\nabla$ such that the Leibnitz rule below holds:
$$
\nabla(f\cdot s)=df\otimes P(s)+f\cdot\nabla(s),
\,\,\,\,\,\,\,\,\,\,\,\forall (f,s)\in \Omega^0(M)\times \Omega^0(\xi).
$$
(See \cite{a2}.)
An {\em ordinary connection} is a pseudoconnection
whose principal homomorphism is the identity.
By classical arguments we shall associate to any pseudoconnection
$\nabla$ an {\em exterior derivative}, that is,
a sequence of linear maps
$d^\nabla:\Omega^k(\xi)\to \Omega^{k+1}(\xi)$
which reduces to $\nabla$ when $k=0$ and satisfies
a Leibnitz rule.
We shall use it to define the
curvature form $F^\nabla:\Omega^0(\xi)\to \Omega^2(\xi)$ as the alternating sum
$$
F^\nabla=P\circ d^\nabla\circ\nabla-d^\nabla\circ P\circ\nabla+d^\nabla\circ\nabla\circ P.
$$
Notice that $F^\nabla$ above reduces to the classical curvature
form if $\nabla$ were an ordinary connection.
We shall prove that $F^\nabla$
is a tensor, that is, $F^\nabla\in Hom(\Omega^0(\xi),\Omega^2(\xi))$,
and explain how to obtain the Abe's curvature \cite{a2} from $F^\nabla$.
We also prove that
$F^\nabla=0$ is a necessary (but not sufficient)
condition for $d^\nabla$ to be a chain complex.
Instead we prove that $F^\nabla=0$ and $d^\nabla\circ d^\nabla\circ\nabla=0$
are necessary and sufficient conditions for $d^\nabla$ to be
a {\em chain $2$-complex}, i.e., $d^\nabla\circ d^\nabla\circ d^\nabla=0$.

\section{Results}

\noindent
Let $\xi,\eta$ be real smooth vector bundles over
a differentible manifold $M$.
An {\em $O$-derivative operator
from $\xi$ to $\eta$ with principal homomorphism $P\in HOM(\xi,\eta)$}
is an $I\!\! R$-linear map $\nabla:\Omega^0(\xi)\to\Omega^0(\eta)$
satisfying the Leibnitz rule
$$
\nabla(f\cdot s)=df\otimes P(s)+f\cdot\nabla(s),
\,\,\,\,\,\,\,\,\,\,\,\forall (f,s)\in \Omega^0(M)\times \Omega^0(\xi).
$$
(See \cite{a1}.)
Notice that a pseudoconnection of $\xi$ is nothing
but an $O$-derivative operator from $\xi$ to itself.

As in \cite{a1} we denote by
$O(\xi,\eta;P)$ the set whose elements are the $O$-derivative operators
with principal homomorphism $P$
from $\xi$ to $\eta$.
We even write $O(\xi;P)$ instead of $O(\xi,\xi;P)$ and
define
$$
O(\xi,\eta)=\bigcup_PO(\xi,\eta;P)
\,\,\,\,\,\,\,\mbox{ and }
\,\,\,\,\,\,\,\,
O(\xi)=\bigcup_PO(\xi;P).
$$
Every $\alpha\in HOM(\xi,\eta)$ induces an alternating product
$$
\wedge_\alpha:\Omega^k(M)\times \Omega^l(\xi)\to\Omega^{k+l}(\eta)
$$
defined at the generators by
$$
\beta\wedge_\alpha(\omega\otimes s)=(\beta\wedge\omega)\otimes\alpha(s).
$$
If $\eta=\xi$ and $\alpha$ is the identity,
then $\wedge_\alpha$ reduced to the ordinary
alternating product $\wedge$ (\cite{mt}).

\begin{lemma}
\label{l1}
For every $\nabla\in O(\xi,\eta;P)$ there is a unique
sequence of linear maps
$d^\nabla:\Omega^k(\xi)\to\Omega^{k+1}(\eta)$, $k\geq 0$,
satisfying the following properties:
\begin{enumerate}
\item
If $k=0$, then
\[
d^\nabla=\nabla.
\]

\item
If $k,l\geq 0$, $\omega\in \Omega^k(M)$ and $S\in \Omega^l(\xi)$, then
\[
d^\nabla(\omega\wedge S)=d\omega\wedge_P S+(-1)^k\omega\wedge d^\nabla S.
\]
\end{enumerate}
\end{lemma}

\noindent
{\bf Proof.}
First define the map $D^\nabla:\Omega^k(M)\times \Omega^0(\xi)\to \Omega^{k+1}(\eta)$ by
$$
D^\nabla(\omega, s)=d\omega\otimes P(s)+(-1)^k\omega\wedge\nabla s,
\,\,\,\,\,\,\,\,\,\forall (\omega,s)\in \Omega^k(M)\times \Omega^0(\xi).
$$
Clearly $D^\nabla$ is
linear and satisfies
$$
D^\nabla(f\cdot \omega,s)=D^\nabla(\omega,f\cdot s)
$$
for all $f\in \Omega^0(M)$ and $(\omega,s)\in \Omega^k(M)\times \Omega^0(\xi)$.
Therefore $D^\nabla$ induces a linear map
$d^\nabla:\Omega^k(\xi)\to \Omega^{k+1}(\eta)$
whose value at the generator $\omega\otimes s$ of $\Omega^k(\xi)$ is given by
$$
d^\nabla(\omega\otimes s)=d\omega\otimes P(s)+(-1)^k\omega\wedge\nabla s.
$$
It follows that $d^\nabla$ and $\nabla$ coincide
at the generators (for $k=0$) by the Leibnitz rule of $\nabla$.
Therefore (1) holds.
The proof of (2) follows as in \cite{mt}.
This ends the proof.

\vspace{5pt}

The sequence $d^\nabla$ in the lemma above will be refered to as
the {\em exterior derivative} of $\nabla\in O(\xi,\eta)$.
Next we state the following definition.

\begin{defi}
Let $\nabla$ be a pseudoconnection with principal
homomorphism $P$ on a vector bundle $\xi$.
We define the following maps
$$
E^\nabla,F^\nabla\colon\Omega^0(\xi)\to \Omega^2(\xi),
\,\,\,\,\,\,
L^\nabla\colon\Omega^0(\xi)\to \Omega^1(\xi)
\,\,\,\,\mbox{ and }
\,\,\,\,\,G^\nabla\colon\Omega^0(\xi)\to\Omega^3(\xi)
$$
by
\begin{itemize}
\item
$E^\nabla=d^\nabla\circ\nabla$;
\item
$F^\nabla=P\circ d^\nabla\circ\nabla- d^\nabla\circ P\circ\nabla
+d^\nabla\circ\nabla\circ P$;
\item
$L^\nabla=P\circ\nabla-\nabla\circ P$;
\item
$G^\nabla=d^\nabla\circ d^\nabla\circ\nabla$.
\end{itemize}
The map $F^\nabla$ will be refered to as the {\em curvature form} of
$\nabla$.
\end{defi}

The maps in the definition above are related by the expressions
\begin{equation}
 \label{EFGL}
F^\nabla=P\circ E^\nabla-d^\nabla\circ L^\nabla,
\,\,\,\,\,\,
G^\nabla=d^\nabla\circ E^\nabla.
\end{equation}
As already observed the curvature form $F^\nabla$
of a pseudoconnection $\nabla$ reduces to the classical curvature form
when $\nabla$ is an ordinary connection \cite{mt}.

\begin{thm}
\label{th1}
If $\nabla$ is a pseudoconnection of $\xi$,
then $F^\nabla\in Hom(\Omega^0(\xi),\Omega^2(\xi))$
and $L^\nabla\in Hom(\Omega^0(\xi),\Omega^1(\xi))$.
\end{thm}

\noindent
{\bf Proof.}
It is not difficult to see that $L^\nabla\in Hom(\Omega^0(\xi),\Omega^1(\xi))$.
On the other hand,
$$
P(\omega\wedge S)=\omega\wedge_P S,
\,\,\,\,\,\,\,\forall \omega\in \Omega^1(M),
\forall S\in \Omega^1(\xi)
$$
so
$$
E^\nabla(f\cdot s)=df\wedge L^\nabla(s)+f\cdot E^\nabla(s),
\,\,\,\,\,\,\,\forall f\in \Omega^0(M),
\forall s\in \Omega^0(\xi).
$$
Then, (\ref{EFGL}) implies
$$
F^\nabla(f\cdot s) =
P(E^\nabla(f\cdot s)))-d^\nabla(L^\nabla(f\cdot s))=
P(df\wedge L^\nabla(s)+f\cdot E^\nabla(s))-d^\nabla(f\cdot L^\nabla(s))
$$
$$
=
df\wedge_PL^\nabla(s)-
df\wedge_PL^\nabla(s)
+f\cdot F^\nabla(s)
=f\cdot F^\nabla(s)
$$
$\forall f\in \Omega^0(M)$, $\forall s\in \Omega^0(\xi)$.
Therefore $F^\nabla\in Hom(\Omega^0(\xi),\Omega^2(\xi))$
and we are done.
This ends the proof.

\vspace{5pt}

\begin{lemma}
\label{impo}
If $\nabla$ is a pseudoconnection of a vector bundle $\xi$ and
$i\geq 0$, then
\begin{equation}
\label{equs}
d^\nabla\circ d^\nabla\circ d^\nabla(\omega\otimes s)
=
d\omega\wedge F^\nabla(s)+(-1)^i\omega\wedge G^\nabla(s)
\end{equation}
for every generator $\omega\otimes s\in\Omega^i(\xi)$.
\end{lemma}

\noindent
{\bf Proof.}
First notice that
$$
d^\nabla\circ d^\nabla(\omega\otimes s) =
d^\nabla(d^\nabla(\omega\otimes s))=
d^\nabla(d\omega\otimes P(s)+(-1)^i\omega\wedge\nabla s)=
$$
$$
=
d^2\omega\otimes P^2(s)+(-1)^{i+1}d\omega\wedge\nabla P(s)
+(-1)^i(dw\wedge_P\nabla s+
$$
$$
\qquad +(-1)^i \omega \wedge d^\nabla(\nabla s))
=
(-1)^i\big[d\omega\wedge(P\nabla s-\nabla Ps)+(-1)^i\omega\wedge
d^\nabla(\nabla s)\big].
$$
Therefore
$$
d^\nabla\circ d^\nabla(\omega\otimes s)=
\omega\wedge E^\nabla(s)+(-1)^id\omega\wedge L^\nabla(s).
$$
Applying $d^\nabla$ to this expression we get (\ref{equs}).
The proof follows.

\vspace{5pt}

As is well known the curvature form $F^\nabla$
of an ordinary connection $\nabla$ measures how
the exterior derivative $d^\nabla$ of $\nabla$ deviates from being a
{\em chain complex}, i.e., $d^\nabla\circ d^\nabla=0$.
Indeed, $d^\nabla$ is a chain complex if and only if
$F^\nabla=0$.
However,
the analogous result for pseudoconnections is false
in general by Proposition \ref{counterexample} below.
Despite we shall obtain a pseudoconnection
version of this result based on the following definition.

\begin{defi}
A pseudoconnection $\nabla$ is called:
\begin{enumerate}
\item
{\em strongly flat} if $E^\nabla=0$ and $L^\nabla=0$,
\item
{\em weakly flat} if $F^\nabla=0$ and $G^\nabla=0$.
\end{enumerate}
\end{defi}

For ordinary connections
one has $F^\nabla=E^\nabla$, $L^\nabla=0$ and then
the notions of flatness above coincide
with the classical flatness  \cite{mt}.
The exterior derivative $d^\nabla$ of a pseudoconnection $\nabla$
is said to be a {\em chain $2$-complex}
if $d^\nabla\circ d^\nabla\circ d^\nabla=0$. 
With these definitions we have the following result.

\begin{thm}
\label{th3}
A pseudoconnection $\nabla$ is weakly flat (resp. strongly flat)
if and only if $d^\nabla$ is a chain $2$-complex (resp. chain complex).
\end{thm}

\noindent
{\bf Proof.}
We only prove the result for weakly flat
since the proof for strongly flat is analogous.

Fix a pseudoconnection $\nabla$
with principal homomorphism $P$ on a vector bundle $\xi$.
If $\nabla$ is weakly flat then
$d^\nabla$ is a chain $2$-complex by (\ref{equs}) in Lemma \ref{impo}.
Conversely, if $d^\nabla$ is a chain $2$-complex, then both $d^\nabla\circ d^\nabla\circ d^\nabla$
and $d^\nabla\circ d^\nabla\circ \nabla$ vanish hence
$\omega\wedge F^\nabla(s)=0$
for all exact form $\omega$ of $M$ and all $s\in \Omega^0(\xi)$
by (\ref{equs}) in Lemma \ref{impo}.
Since every form in $M$ is locally a $\Omega^0(M)$-linear combination of alternating product
of exact forms we obtain that $\omega \wedge F^\nabla(s)=0$
for all $k$-form $\omega$ of $M$ ($k\geq 1$) and all $s\in \Omega^0(\xi)$.
From this we obtain that $F^\nabla=0$ so $\nabla$ is weakly flat.
The proof follows.

\vspace{5pt}

The following is a direct corollary of the above theorem.

\begin{clly}
If the exterior derivative $d^\nabla$ of a pseudoconnection $\nabla$ is
a chain complex, then $F^\nabla=0$.
\end{clly}

The converse of the above corollary is false by the following
proposition.

\begin{prop}
\label{counterexample}
There is a pseudoconnection $\nabla$ with $F^\nabla=0$
such that $d^\nabla$ is not a chain complex.
\end{prop}

\noindent
{\bf Proof.}
Choose a suitable vector bundle $\xi$ over $M=I\!\! R^3$,
$\Phi_2,\Phi_3\in HOM(\xi,\xi)$ such that
$\Phi_3\circ\Phi_2\neq\Phi_2\circ\Phi_3$ and three
$1$-forms $\omega_1,\omega_2,\omega_3\in \Omega^1(M)$
such that $\omega_1\wedge\omega_2\wedge\omega_3$ never vanishes.
Define the map
$\nabla:\Omega^0(\xi)\to\Omega^1(\xi)$
by
$$
\nabla s=\omega_1\otimes s+\omega_2\otimes \Phi_2(s)+\omega_3\otimes\Phi_3(s).
$$
We have that $\nabla\in Hom(\Omega^0(\xi),\Omega^1(\xi))$ therefore
$\nabla$ is a pseudoconnection with zero principal homomorphism
so $F^\nabla=0$.
On the other hand,
an straightforward computation yields
$$
G^\nabla(s)=
(\omega_1\wedge\omega_2\wedge\omega_3)\otimes
(\Phi_3\circ\Phi_2-\Phi_2\circ\Phi_3)(s),
\,\,\,\,\,\,\,\,\forall s\in \Omega^0(\xi)
$$
therefore $G^\nabla\neq0$ and so $\nabla$
is not weakly flat. Then, $d^\nabla$ cannot be a chain complex
by Theorem \ref{th3} since a chain complex is necessarily a chain $2$-complex.
This ends the proof.

\vspace{5pt}

To finish we explain how
the Abe's curvature \cite{a2} can be obtained from
the curvature form $F^\nabla$.
For this we need some short definitions (see \cite{mt}).

Given a vector bundle $\xi$ over $M$ and $k$ vector  fields
$X_1,\cdots ,X_k\in \chi^\infty(M)$
we define the evaluation map $Ev_{X_1,\cdots, X_k}:\Omega^k(\xi)\to \Omega^0(\xi)$
by defining
$$
Ev_{X_1,\cdots, X_k}(\omega\otimes s)=
w(X_1,\cdots,X_k)\cdot s
$$
at each generator $\omega\otimes s\in \Omega^k(\xi)$.
If $\nabla\in O(\xi)$ and $X,Y\in \chi^\infty(M)$ we define
$\nabla_X:\Omega^0(\xi)\to\Omega^0(\xi)$
by $$
\nabla_Xs=Ev_X(\nabla s)
$$
and
$F^\nabla_{X,Y}:\Omega^0(\xi)\to \Omega^0(\xi)$
by
$$
F^\nabla_{X,Y}=
Ev_{X,Y}(F^\nabla(s)),
\,\,\,\,\,\,\forall s\in \Omega^0(\xi).
$$
\begin{thm}
\label{th2}
If $\nabla\in O(\xi;P)$ then
$$
F^\nabla_{X,Y}(s) =  \nabla_X\nabla_Y(Ps)-\nabla_Y\nabla_X(Ps)
- \nabla_XP(\nabla_Ys)+P\nabla_X\nabla_Ys+
$$
$$
+\nabla_YP(\nabla_Xs)-
P\nabla_Y\nabla_Xs-
P\left(\nabla_{[X,Y]}P(s)\right),
$$
for all $X,Y\in \chi^\infty(M)$ and all $s\in \Omega^0(\xi)$.
\end{thm}

\noindent
{\bf Proof.}
By definition we have
\begin{equation}
\label{eq1}
F^\nabla_{X,Y}(s)=
\end{equation}
$$
Ev_{X,Y}(P(d^\nabla(\nabla s)))-Ev_{X,Y}(d^\nabla(P(\nabla s)))
+Ev_{X,Y}(d^\nabla(\nabla(Ps))).
$$

Let us compute the three sumands
separated way.
First of all if $\omega\otimes s\in \Omega^1(\xi)$ is a generator then
\begin{equation}
\label{eq2}
Ev_{X,Y}(d^\nabla(\omega\otimes s))=
dw(X,Y)\cdot P(s)-\omega(X)\cdot\nabla_Ys+\omega(Y)\cdot\nabla_Xs.
\end{equation}

Now, as $\nabla s\in \Omega^1(\xi)$ and 
$\{\omega\otimes s':(\omega,
s')\in \Omega^1(M)\times \Omega^0(\xi)\}$ is a generating set of
$\Omega^1(\xi)$ we obtain
\begin{equation}
\label{eq3}
\nabla s=\sum_{r=1}^k\omega_r\otimes s_r,
\end{equation}
for some $(\omega_r, s_r) \in \Omega^1(M)\times \Omega^0(\xi)$,
$r=1,\cdots ,k$.
Then (\ref{eq2}) implies
\begin{equation}
\label{eq4}
Ev_{X,Y}(d^\nabla(\nabla s))=
\end{equation}
$$
\sum_{r=1}^k\big\{dw_r(X,Y)\cdot P(s_r)-\omega_r(X)\cdot\nabla_Ys_r+\omega_r(Y)\cdot\nabla_Xs_r\big\}.
$$
On the other hand, (\ref{eq3}) yields
$$
\nabla_Xs=\sum_{r=1}^k\omega_r(X)\cdot s_r
$$
therefore
$$
\nabla_Y\nabla_Xs=
\sum_{r=1}^k
\left\{d[\omega_r(X)](Y)\cdot P(s_r)+
\omega_r(Y)\cdot \omega_r(X)\cdot\nabla_Ys_s\right\}.
$$
But
$\nabla_{[X,Y]}s=\displaystyle \sum_{r=1}^k\omega_r([X,Y])\cdot s_r$,
so
$$
P\left(\nabla_{[X,Y]}s\right)=\sum_{r=1}^k\omega_r([X,Y])\cdot P(s_r)
$$
and then
\begin{equation}
\label{eq5}
Ev_{X,Y}(d^\nabla(\nabla s))=
\nabla_X\nabla_Y s-\nabla_Y\nabla_X s-P\left(\nabla_{[X,Y]}s\right)
\end{equation}
because of (\ref{eq4}).
Replacing $s$ by $P(s)$ in (\ref{eq5}) we obtain
\begin{equation}
\label{eq6}
Ev_{X,Y}(d^\nabla(\nabla P(s)))=
\nabla_X\nabla_YP(s)-\nabla_Y\nabla_XP(s)-
P\left(\nabla_{[X,Y]}P(s)\right).
\end{equation}
Besides (\ref{eq3}) implies
$$
P(\nabla s)=
\sum_{r=1}^k\omega_r\otimes P(s_r)
$$
thus
$$
Ev_{X,Y}(d^\nabla(P\nabla s))=
\sum_{r=1}^k Ev_{X,Y}(d^\nabla(\omega_r\otimes P(s_r)))
=
$$
$$
\sum_{r=1}^k
\big[dw_r(X,Y)\cdot P^2(s_r)-
\omega_r(X)\nabla_XP(s_r)
+\omega_r(Y)\cdot \nabla_XP(s_r)\big]
$$
and then
$$
\nabla_YP(\nabla_Xs)=
\sum_{r=1}^k\left\{d[w_r(X)](Y)\cdot P^2(s_r)+w_r(X)\cdot\nabla_YP(s_r)\right\}.
$$
As
$P^2\left(\nabla_{[X,Y]}s\right)=
\displaystyle \sum_{r=1}^k\omega_r([X,Y])\cdot P^2(s_r)$
we obtain
\begin{equation}
\label{eq7}
Ev_{X,Y}(d^\nabla(P(\nabla s)))=
\nabla_XP(\nabla_Ys)-\nabla_YP (\nabla_Xs)-
P^2\left(\nabla_{[X,Y]}s\right).
\end{equation}
As the maps $P$ and $Ev_{X,Y}$ commute
we can apply
$P$ to (\ref{eq5}) and use (\ref{eq1}), (\ref{eq6}), (\ref{eq7})
to obtain the result.

\vspace{5pt}

\begin{rk}
$F^\nabla_{X,Y}(s)$ in Theorem \ref{th2} is
the curvature $K(\nabla)_{X,Y}(s)$ defined in \cite{a2} p. 328.
\end{rk}

\section*{Acknowledgements}

\noindent
The first author was partially supported by CNPq, FAPERJ and PRONEX/DS from Brazil.

\end{document}